\documentclass[12pt,leqno]{article}
\usepackage{amsfonts}
\usepackage{amsmath, amsthm, amsfonts, amssymb, color}
\usepackage{mathrsfs}
\usepackage{color}
\newtheorem{thm}{Theorem}[section]

\newtheorem{exa}[thm]{Example}

\theoremstyle{definition}

\newcommand{\scr}[1]{\mathscr #1}
\definecolor{wco}{rgb}{0.5,0.2,0.3}

\numberwithin{equation}{section} \theoremstyle{remark}

\newcommand{\ua}{\uparrow}

\title{{\bf Perturbations of Functional Inequalities for L\'evy Type Dirichlet Forms}\footnote{Supported in
 part by  Lab. Math. Com. Sys., NNSFC(11131003 and 11201073), the Program for New Century Excellent Talents in Universities
of Fujian (No.\ JA12053), the Program for Nonlinear Analysis and Its Applications (No.\ IRTL1206), and the Portuguese Science Foundation (FCT)
for the project ``Probabilistic approach to finite and infinite
dimensional dynamical systems'' (No.\ PTDC/MAT/104173/2008).} }
\author{
{\bf   Xin Chen$^{b)}$, Feng-Yu Wang$^{a), c)}$  and Jian Wang$^{d)}$}\\
\footnotesize{$^{a)}$School of Mathematical Sciences,
Beijing Normal
University, Beijing 100875, China}\\
\footnotesize{$^{b)}$Grupo de Fisica Matematica, Universidade de
Lisboa, Av Prof Gama Pinto 2, Lisbon 1649-003, Portugal}\\
 \footnotesize{$^{c)}$Department of Mathematics,
Swansea University, Singleton Park, SA2 8PP, United Kingdom}\\
\footnotesize{$^{d)}$School of Mathematics and Computer Science,
Fujian Normal University, Fuzhou 350007, China}\\
\footnotesize{chenxin\_217@hotmail.com, wangfy@bnu.edu.cn,
F.-Y.Wang@swansea.ac.uk,
 jianwang@fjnu.edu.cn}}
\begin{document}
\allowdisplaybreaks
\def\R{\mathbb R}  \def\ff{\frac} \def\ss{\sqrt} \def\B{\mathbf
B}
\def\N{\mathbb N} \def\kk{\kappa} \def\m{{\bf m}}
\def\ee{\varepsilon}
\def\dd{\delta} \def\DD{\Delta} \def\vv{\varepsilon} \def\rr{\rho}
\def\<{\langle} \def\>{\rangle} \def\GG{\Gamma} \def\gg{\gamma}
  \def\nn{\nabla} \def\pp{\partial} \def\E{\scr E}
\def\d{\text{\rm{d}}} \def\bb{\beta} \def\aa{\alpha} \def\D{\scr D}
  \def\si{\sigma} \def\ess{\text{\rm{ess}}}
\def\beg{\begin} \def\beq{\begin{equation}}  \def\F{\scr F}
\def\Ric{\text{\rm{Ric}}} \def\Hess{\text{\rm{Hess}}}
\def\e{\text{\rm{e}}} \def\ua{\underline a} \def\OO{\Omega}  \def\oo{\omega}
 \def\tt{\tilde} \def\Ric{\text{\rm{Ric}}}
\def\cut{\text{\rm{cut}}} \def\P{\mathbb P} \def\ifn{I_n(f^{\bigotimes n})}
\def\C{\scr C}      \def\aaa{\mathbf{r}}     \def\r{r}
\def\gap{\text{\rm{gap}}} \def\prr{\pi_{{\bf m},\varrho}}  \def\r{\mathbf r}
\def\Z{\mathbb Z} \def\vrr{\varrho} \def\ll{\lambda}
\def\L{\scr L}\def\Tt{\tt} \def\TT{\tt}\def\II{\mathbb I}
\def\i{{\rm in}}\def\Sect{{\rm Sect}}  \def\H{\mathbb H}
\def\M{\scr M}\def\Q{\mathbb Q} \def\texto{\text{o}} \def\LL{\Lambda}
\def\Rank{{\rm Rank}} \def\B{\scr B} \def\i{{\rm i}} \def\HR{\hat{\R}^d}
\def\to{\rightarrow}\def\l{\ell}
\def\EE{\scr E}
\def\A{\scr A}
\def\BB{\scr B}

\maketitle

\begin{abstract} Perturbations of super  Poincar\'{e} and weak Poincar\'{e} inequalities for L\'evy type   Dirichlet forms are studied. When the range of jumps is finite our results are natural extensions to the corresponding ones    derived earlier for diffusion processes; and  we show that the study for the situation with infinite range of jumps is essentially different. Some examples are presented to illustrate the optimality of our results.
     \end{abstract} \noindent
 AMS subject Classification:\  60J75, 47G20, 60G52.   \\
\noindent
 Keywords:  Super   Poincar\'{e} inequality;  weak  Poincar\'{e} inequality; L\'evy type Dirichlet form; perturbation.
 \vskip 2cm

\section{Introduction}

Functional inequalities of Dirichlet forms are powerful tools in the study of Markov semigroups and spectral theory of Dirichlet operators, see \cite{B,D,G,Wbook} for accounts on functional inequalities and applications. To establish a functional inequality, one often needs to verify some conditions on the generator, for instances the Bakry-Emery curvature condition in the diffusion setting or the Lyapunov condition in a
 general setting, see e.g. \cite{CGWW,GLWN,Wbook}. Since these conditions exclude generators with less regular coefficients, to establish functional inequalities in a more general setting one treats the singularity part as a perturbation. So, it is important to investigate perturbations of functional inequalities.

It is a classical result that all Poincar\'{e}-Sobolev type
inequalities are stable under bounded perturbations, see e.g.\
\cite[Section 3.3]{Ch} and references therein. Under some regularity
conditions, functional inequalities could also be stable under
unbounded perturbations. For instance, in the diffusion setting
(i.e.\ the underlying Dirichlet form is local) the stability of the
log-Sobolev inequality was proved in \cite[Section 3]{AS}  under an
exponential integrability condition on the gradient of the
perturbation term. As for less regular perturbations, sharp growth
conditions have been presented in \cite{BLW}  for perturbations of
super Poincar\'e and  weak Poincar\'e inequalities. Note that these
two kinds of functional inequalities are general enough to cover all
Poincar\'e/Sobolev/Nash type inequalities, and thus, have a broad
range of applications. However, perturbations of functional
inequalities for non-local Dirichlet forms are much less known in
the literature.

Recently, explicit sufficient conditions were derived in
\cite{CW1, CW2, MRS, WW, WJ13} for functional inequalities of
stable-like Dirichlet forms. The aim of this paper is to extend
perturbation results derived in \cite{BLW} to the non-local setting,
so that combining with the existing sufficient conditions we are
able to establish functional inequalities for more general Dirichlet
forms. Due to the lack of the
chain rule, the present study   is much more
complicated than in the diffusion setting. Nevertheless, we are able to present some relatively
clean and new perturbation results, which are sharp as illustrated
by some examples latter on, and when the range of jumps is finite,
are natural extensions to the corresponding results derived earlier
for diffusion processes.

Let $(E,d)$ be a Polish space equipped with the Borel $\si$-field $\F$ and a probability measure $\mu$. Let $\B(E)$ be the set of all measurable functions on $E$, and let $\B_b(E)$ be the set of all bounded elements in $\B(E).$ Let $q\in\B(E\times E)$ be  non-negative   with $q(x,x)=0, x\in E,$ such that
\begin{equation}\label{intro-02}
\ll:= \sup_{x\in E}\int_{E} \big(1\land d(x,y)^2\big)
q(x,y) \mu(\d y)<\infty.\end{equation}Then
\begin{equation}\label{intro-01}\beg{split}
&\GG(f,g)(x):= \int_{E} (f(x)-f(y))(g(x)-g(y)) q(x,y)\mu(\d y),\ \ \forall x \in E,\\
&f,g\in \scr A:= \big\{f\in \B_b(E):  \GG(f,f)\in \B_b(E)\big\}\end{split}\end{equation} gives rise to a non-negatively definite bilinear map from $\A\times\A$ to $\B_b(E).$  We assume that $\A$ is dense in $L^2(\mu)$. Then it is standard that the form
\beq\label{intro-03} \beg{split}\E(f,g):&= \mu(\GG(f,g))\\
 &=\int_{E\times E} (f(x)-f(y))(g(x)-g(y))q(x,y)\mu(\d y)\mu(\d x),\ \ f,g\in \A\end{split}\end{equation}  is closable in $L^2(\mu)$ and its closure $(\E,\D(\E))$ is a symmetric conservative Dirichlet form,
see e.g.\ \cite[Example 1.2.6]{FOT}. A typical example of the
framework is the $\aa$-stable-like Dirichlet form, where $E=\R^m$
with $d(x,y)=|x-y|$,  and
\begin{equation*}\label{intro-04} q(x,y)= \ff{\tt q (x,y)}{u(y)|x-y|^{m+\aa}},\ \ \mu(\d x)= u(x)\d x\end{equation*} for some $\aa\in (0,2)$,   non-negative $\tt q\in \B_b(\R^m\times\R^m)$, and positive $u\in\B(\R^m)$ such that $\mu(\d x)$ is a probability measure.  In this case we have $\A\supset C_0^2(\R^m)$ which is dense in $L^2(\mu)$ for any probability measure $\mu$ on $\R^m$.

To investigate perturbations of functional inequalities using growth conditions as in \cite{BLW}, we fix a point $o\in E$ and  denote $\rr(x)=d(o,x), x\in E.$ Now, for
 a $\rho$-locally bounded measurable function $V$ on $E$ (i.e. $V$ is bounded on the set $\{x \in E: \ \rho(x)\le r\}$ for all $r>0$)  such that
 $\mu(\e^V)=1,$ let $\mu_V(\d x) =   \e^{V(x)} \mu(\d x)$.
 Since for every $f \in \A$, $\GG(f,f)$ given by (\ref{intro-01}) is a bounded measurable function on $E$, we have
$$
 \EE_V(f,f):=\int_{E}\GG(f,f)(x)\,\mu_V(\d x) <\infty.
$$
Again by the argument in \cite[Example 1.2.6]{FOT},  the form
$$  \E_V(f,g):= \mu_V(\GG(f,g))=\int_{E\times E} (f(x)-f(y))(g(x)-g(y))q(x,y)\,\mu(\d y)\,\mu_V(\d x)$$ defined for $f,g\in \A$ is closable in $L^2(\mu_V)$ and its closure $(\E_V,\D(\E_V))$ is a symmetric conservative Dirichlet form.

In this paper, we shall assume that $\EE$ satisfies a functional
inequality and then search for conditions on $V$ such that $\EE_V$
satisfies the same type of functional inequality. We now briefly
explain the motivation for the study of the perturbed Dirichlet form
${\EE_V}$, which only changes the underlying measure $\mu$ and keeps
the non-negatively definite bilinear map $\GG(f,f)$ unchanged. Thus, the perturbation occurs essentially in the drift part
of the generator of the associated Markov process. Indeed, take $E=\R^m$ and $d(x,y)=|x-y|$ for example. Let $(\mathcal L,C_0^\infty(\R^m))$ be the pre-generator of a Markov process on $\R^m$ with
invariant probability measure $\mu$ and square field
$$\GG(f,f)= \ff 1 2 \mathcal L f^2- f\mathcal Lf,\ \ f\in C_0^\infty(\R^m).$$ Consider the operator $\mathcal L_b:=\mathcal L+ b\cdot\nn$ for some measurable
$b: \R^m\to\R^m$. Then the associated square field of $\mathcal L_b$
is still $\GG(f,f)$, but the invariant probability measure (if
exists) might be different from $\mu$, which is then regarded as a
perturbed measure of $\mu$.

In the following two sections, we will investigate  perturbations
of the super Poincar\'e inequality and the weak Poincar\'e
inequality respectively. As mentioned above that known results on
functional inequalities for the non-local Dirichlet form $\EE_V$
(see \cite{CW2, MRS,WW}) require regularity conditions on $V,$  and
in our results we only use growth conditions of $V.$ So, this paper does
provide new results, which lead to examples of  non-local
Dirichlet forms with singular potentials  satisfying functional
inequalities. For simplicity, throughout the paper we will denote
$\GG(f)=\GG(f,f),$ $\EE(f)=\EE(f,f)$ and $\EE_V(f)=\EE_V(f,f).$

\section{Perturbations of the Super Poincar\'e inequality}

We consider the following super Poincar\'e inequality introduced in
\cite{W00a, W00b}:
 \beq\label{sp} \mu(f^2)\le r \EE(f) +\bb(r)\mu(|f|)^2,\ \ r>0, f\in \scr{D}(\EE),\end{equation}
 where $\bb: (0,\infty)\to (0,\infty)$ is a decreasing function.
 Note that if \eqref{sp} holds for some non-decreasing function $\bb$,
 it holds also for the decreasing function
 $\tt\bb(r):=\inf_{s\in (0,r]} \bb(s)$ in place of $\bb(r).$

In the first two subsections, we consider Dirichlet forms with jumps of finite range and of infinite range respectively. Then the optimality of  results derived in these two situations is   illustrated in \S 2.3
by specific examples. Finally, to compare with perturbation results derived in the diffusion setting, a variant condition is introduced in \S 2.4.

\subsection{Jumps of finite range}

To establish a super Poincar\'e inequality for $\EE_V$, we  need the
following   quantities. For any $n\ge 1$ and $k\ge 1$, let
\beg{equation*}\beg{split} &K_{n,k}(V)= \sup_{\rho(x)\le n+2}V(x)-\inf_{\rho(x)\le n+k+2}V(x),\\
 &J_{n,k}(V)= \sup_{\rho(x)\le n+1}V(x)-2\inf_{\rho(x)\le n+k+2}V(x),\\
 &\vv_{n,k}(V) = \sup_{m\ge n} \Big\{\bb^{-1}\big(1/[2\mu(\rho>m-1)]\big) \e^{K_{m,k}(V)}\Big\},
 \end{split}\end{equation*}
 where
 $\bb^{-1}(s):=\inf\{r>0:\bb(r)\le s\}$ for $s>0$,  with $\inf\emptyset :=\infty$ by convention.
 Roughly speaking, the parameter $\vv_{n,k}(V)$ is a factor to measure the oscillation condition on $V$ and the decay condition on the rate function $\bb$ in the given super Poincar\'{e} inequality $(\ref{sp})$.

When the jump is of finite range, i.e. there exists $k_0>0$ such
that $q(x,y)=0$ for $d(x,y)> k_0$, we have the following   result corresponding to
  \cite[Theorem 3.1]{BLW}   in the diffusion setting.

\beg{thm} \label{finite-sp} Assume that $(\ref{sp})$ holds and there
exists $k_0\ge 1 $ such that $q(x,y)=0$ for  $d(x,y)> k_0$.
\begin{itemize}
\item[$(1)$] If $\inf\limits_{n\ge1,k\ge k_0}\vv_{n,k}(V)=0,$ then the super Poincar\'{e} inequality
\begin{equation}\label{th-sp-1}\mu_V(f^2)\le
r\EE_V(f)+\bb_V(r) \mu_V(|f|)^2,\ \ r>0,\ \ f\in
\scr{D}(\EE_V)\end{equation}holds with
\beg{equation*}\beg{split} \bb_V(r):=\inf\Big\{&(1+8\lambda r')\e^{J_{n,k}(V)}\bb(s): s>0, r'\in (0,r], n\ge 1, k\ge k_0\\
&\quad {\rm such\ that}\ 8\vv_{n,k}(V)+s\e^{K_{n,k}(V)}\le
\frac{r'}{2+16\ll r'} \Big\}.\end{split}\end{equation*}
\item[$(2)$] If $\inf\limits_{n\ge 1,k\ge k_0}\vv_{n,k}(V)<\infty,$ then
the Poincar\'{e} inequality \beq\label{th-sp-2-1}\mu_V(f^2)\le
C\EE_V(f) +\mu_V(f)^2,\ \ f\in\D(\EE_V)\end{equation} holds for some
constant $C>0$.
\end{itemize}
\end{thm}

To prove this result, we shall adopt a split argument as in \cite[Theorem 3.1]{BLW}  by
estimating $\mu(f^21_{\{\rr\ge n\}})$ and $\mu(f^21_{\{\rr\le n\}})$
respectively. For later use in  the proof of Theorem $\ref{th-sp}$,  the following two lemmas  allows the jump to have  infinite range.
Unlike in the local setting where the chain rule is
available, in the present situation  the uniform norm $\|f\|_\infty$
will appear in our estimates.
Below we simply denote $K_{n,k}=K_{n,k}(V), J_{n,k}=J_{n,k}(V),
Z_n=Z_n(V),\vv_{n,k}=\vv_{n,k}(V), \zeta_n=\zeta_n(V)$ and
$t_{i,n,k}=t_{i,n,k}(\dd,V).$ Moreover, let
\beg{equation*}\beg{split} &\gg_{n,k}  =\int_{\{d(x,y)>k, \rr(y)\ge n-1\}}q(x,y)\mu(\d y)\mu(\d x),\\
&\eta_{n,k} =\int_{\{\rr(x)>n+k+2, \rr(y)\le n+1\}}q(x,y)\mu(\d
y)\mu(\d x).\end{split}\end{equation*} By  \eqref{intro-02}, we have
  $\gg_{n,k}+\eta_{n,k} \downarrow 0$ as $n\uparrow \infty$ for all
 $k\ge 1.$

\beg{lem}\label{L3.1} For any $n\ge 1$, $k\ge 1$ and $f\in \A$,
$$\mu_V(f^21_{\{\rr\ge n\}})\le 12 \vv_{n,k} \EE_V(f)+ 128\ll\vv_{n,k} \mu_V(f^2)+96 \zeta_{n} \gg_{n,k} \|f\|_\infty^2.$$\end{lem}

\beg{proof} If $f|_{\{\rho\le n-1\}}=0$,
  by the Cauchy-Schwarz inequality we have
\begin{equation*}
\mu(|f|)^2=\mu\big(|f| 1_{\{\rho>n-1\}}\big)^2\le
\mu(\rho>n-1)\mu(f^2).
\end{equation*}
Substituting this into  \eqref{sp}  with $r=
\bb^{-1}\big(\frac{1}{2\mu(\rho>n-1)}\big)$  we obtain
\beq\label{proof-e-01} \mu(f^2)\le 2\bb^{-1}
\big(1/[2\mu(\rho>n-1)]\big)\EE(f),\quad f\in \A, \ f|_{\{\rho\le
n-1\}}=0.\end{equation}

To apply (\ref{proof-e-01}) for general $f\in\A$, we consider
$fl_{n}$ instead of $f$, where $l_{n}:= h_{n}\circ\rr$ for some
function  $h_n\in C^\infty([0,\infty); [0,\infty))$   satisfying
\begin{equation*} h_{n}(s)
\begin{cases}
& =1,\ \ \ \ \ \ \ \quad \,\, n\le s\le n+1,\\
& \in [0,1],\ \ \ \ n-1\le s<n \,\,\,\,{\textrm{ or}}\,\, \,\,n+1<s\le n+2,\\
& =0,\ \ \ \ \ \ \ \ \quad s<n-1\,\,\,\,{\textrm {or}}\,\,\,\,s>n+2,
\end{cases}
\end{equation*}
and $$\sup_{s\in[0,\infty)}|h'_{n}(s)|\le 2.$$    It is easy to see
that
$$|l_n(x)-l_n(y)|\le 2(1\wedge d(x,y))\big\{1_{\{\rr(x)\in (n-2,n+3)\}}+1_{\{\rr(y)\in (n-1,n+2), \rr(x)\notin (n-2,n+3)\}}\big\}.$$ This implies
\beg{equation}\label{proof-e-02}\beg{split} &\GG(fl_n)(x) \\
& \le 2\int_{E} \big\{l_n(y)^2 (f(x)-f(y))^2
 +  f(x)^2   (l_n(x)-l_n(y))^2\big\}q(x,y)\mu(\d y)\\
&\le 2 \int_{\{ \rho(y)\in (n-1,n+2)\}} (f(x)-f(y))^2 q(x,y)\mu(\d y)+8\ll f(x)^21_{\{\rr(x)\in (n-2,n+3)\}}\\
  &\quad + 8f^2(x) 1_{\{\rr(x)\notin(n-2,n+3)\}}\int_{\{\rr(y)\in (n-1,n+2)\}}   q(x,y)\mu(\d y). \end{split}\end{equation}
Since $$1_{\{\rr(x)\notin(n-2,n+3)\}}\int_{\{\rr(y)\in (n-1,n+2)\}}
q(x,y)\mu(\d y)\le 1_{\{\rr(x)\notin(n-2,n+3)\}}\int_{\{d(x,y)>1\}}
q(x,y)\mu(\d y)\le \lambda,$$ we have $\GG(fl_n)\in \B_b(E)$, and
$fl_n \in \A$.

Let
  $$\dd_{n,k} =  \e^{K_{n,k}}  \bb^{-1} \big( 1/[2\mu(\rho>n-1)]\big),\
  \theta_{n}:=\e^{\sup_{\rho\le n+2}V}  \bb^{-1} \big( 1/[2\mu(\rho>n-1)]\big).$$
 Combining
(\ref{proof-e-01}) with \eqref{proof-e-02}, and noting that
$\rr(y)\in (n-1,n+2)$ and $\rr(x)\ge n+k+2$ imply $d(x,y)>k$, and
\begin{equation*}
\begin{split}
&\int_{\{\rr(x)>n+k+2, \rho(y)\in (n-1,n+2)\}} (f(x)-f(y))^2 q(x,y)\mu(\d y)\mu(\d x)\\
&\le 4 \|f\|_\infty^2 \int_{\{\rr(x)>n+k+2, \rr(y) \in (n-1,n+2),d(x,y)>k\}}q(x,y)\mu(\d y)\mu(\d x),\\
\end{split}
\end{equation*}
 we obtain
\beq\label{proof-e-03}\beg{split}&\mu_V(f^2l_n^2) \le
\e^{\sup_{\rho\le n+2}V}\mu(f^2l_n^2)\\
&\le 2\e^{\sup_{\rho\le n+2}V}  \bb^{-1} \big( 1/[2\mu(\rho>n-1)]\big)\EE(fl_n)\\
&\le 2\e^{\sup_{\rho\le n+2}V}  \bb^{-1} \big( 1/[2\mu(\rho>n-1)]\big)\mu(1_{\{\rho\le n+k+2\}}\GG(fl_n))\\
&\quad + 2\e^{\sup_{\rho\le n+2}V}  \bb^{-1} \big( 1/[2\mu(\rho>n-1)]\big)\mu(1_{\{\rho> n+k+2\}}\GG(fl_n))\\
  &\le   4 \dd_{n,k}\int_{\{\rr(x)\le n+k+2, \rr(y)\in (n-1,n+2)\}}(f(x)-f(y))^2q(x,y)\mu(\d y)\mu_V(\d x) \\
 &\quad +16\dd_{n,k}\ll \mu_V\big(1_{\{\rr\in (n-2,n+3)\}}f^2\big)\\
 &\quad+16\dd_{n,k} \int_{\{\rr(x)\le n+k+2, \rr(x) \notin (n-2,n+3),\rr(y)\in (n-1,n+2)\}}f(x)^2q(x,y)\mu(\d y)\mu_V(\d x)\\
 &\quad+32 \theta_{n}\|f\|_\infty^2\int_{\{\rr(y)\in (n-1,n+2),d(x,y)>k\}}q(x,y)\mu(\d y)\mu(\d x). \end{split}\end{equation}
Noting that $\vv_{n,k}=\sup_{m\ge n} \dd_{m,k}$, $\zeta_n=\sup_{m
\ge n} \theta_m$, $\sum_{m=1}^\infty 1_{\{\rr(y)\in (m-1,m+2)\}}\le
3$ and $$1_{\{\rr(x)\le n+k+2, \rr(x) \notin (n-2,n+3),\rr(y)\in
(n-1,n+2)\}} \le 1_{\{d(x,y)>1, \rr(y)\in (n-1,n+2)\}},$$ and taking
summations in (\ref{proof-e-03}) from $n$, we arrive at
 \beg{equation*}\beg{split}  & \mu_V(f^2 1_{\{\rho\ge n\}}) \le
\sum_{m=n}^\infty \mu_V(f^2l_m^2)\\
&\le 12 \vv_{n,k}\int_{E\times E}(f(x)-f(y))^2q(x,y)\mu(\d y)\mu_V(\d x)   +80\vv_{n,k}\ll \mu_V(f^2)\\
 &\quad+48\vv_{n,k}\int_{\{d(x,y)>1\}}f(x)^2q(x,y)\mu(\d y)\mu_V(\d x) \\
 &\quad+
 96\zeta_n\|f\|_\infty^2\int_{\{\rr(y)\ge n-1,d(x,y)>k\}}q(x,y)\mu(\d y)\mu(\d x)\\
&\le 12 \vv_{n,k} \EE_V(f)+ 128\ll\vv_{n,k} \mu_V(f^2)+96 \zeta_{n}
\gg_{n,k} \|f\|_\infty^2.\end{split}\end{equation*}
\end{proof}

\beg{lem} \label{L3.2} For any $n,k\ge 1$, $s>0$ and $f\in \A$,
$$\mu_V(f^21_{\{\rr\le n\}}) \le  2s \e^{K_{n,k}}\EE_V(f) +16\ll s\e^{K_{n,k}}\mu_V(f^2) + 16 s \e^{Z_{n}}\eta_{n,k} \|f\|_\infty^2
 +\bb(s) \e^{J_{n,k}} \mu_V(|f|)^2.$$\end{lem}

\beg{proof} Let $\phi_n:[0,\infty)\to [0,\infty)$ be a smooth
function such that
\begin{equation*} \phi_n(s) \begin{cases}
& =1,\ \ \ \ \ \ \ \quad \,\, s\le n,\\
& \in [0,1],\ \ \ \ n\le s<n+1,\\
& =0,\ \ \ \ \ \ \ \ \quad s>n+1,
\end{cases}
\end{equation*}
and $$\sup_{s \in [0,\infty)}|\phi'_{n}(s)|\le 2.$$ Set $g_n
=\phi_n\circ\rho.$ Then $g_n\in\A$ and
$$|g_n(x)-g_n(y)|\le   2(1\wedge d(x,y))\big(1_{\{\rr(x)\le n+1\}}+1_{\{\rr(x)>n+1,\rr(y)\le n+1\}}\big).$$
So, similarly to (\ref{proof-e-02}) we have
\beg{equation*}\beg{split} \GG(fg_n)(x) &\le 2 \int_{\{\rr(y)\le n+1\}}(f(x)-f(y))^2q(x,y)\mu(\d y) \\
&\quad + 8\ll f^2(x)1_{\{\rr(x)\le n+k+2\}} + 8 f^2(x)
1_{\{\rr(x)>n+k+2\}} \int_{\{\rr(y)\le n+1\}} q(x,y) \mu(\d y).
\end{split}\end{equation*}
Note that $\rr(x)> n+k+2$ and $\rr(y)\le n+1$ imply that $d(x,y)>k$,
%\begin{equation*}
%\begin{split}
%&\int_{\{n+2<\rr(x)\le n+k+2,\rr(y)\le n+1\}} f^2(x)q(x,y) \mu(\d y)\mu_V(\d x)\le \lambda \mu_V(f^2),
%\end{split}
%\end{equation*}
and
\begin{equation*}
\begin{split}
&\int_{\{\rr(x)>n+k+2,\rho(y)\le n+1\}} (f(x)-f(y))^2 q(x,y)\mu(\d y)\mu(\d x)\\
&\le 4 \|f\|_\infty^2 \int_{\{\rr(x)>n+k+2, \rr(y)\le
n+1,d(x,y)>k\}}q(x,y)\mu(\d y)\mu(\d x).
\end{split}
\end{equation*}
Combining all the estimates above with \eqref{sp},   we obtain
\beg{equation*}\beg{split} &\mu_V(f^21_{\{\rho\le n\}})\\
&  \le \mu_V(f^2g_n^2)
 \le \e^{\sup_{\rho\le n+1}V}\mu(f^2g_n^2)\\
&\le \e^{\sup_{\rho\le n+1}V}\Big\{s\mu(\GG(fg_n)) +\bb(s) \mu(|fg_n|)^2\Big\}\\
&= \e^{\sup_{\rho\le n+1}V}\Big\{s\mu(1_{\{\rho\le
n+k+2\}}\GG(fg_n))
 +s\mu(1_{\{\rho> n+k+2\}}\GG(fg_n))+\bb(s) \mu(|fg_n|)^2\Big\}\\
&\le s\e^{K_{n,k}} \mu_V(1_{\{\rho\le n+k+2\}}\GG(fg_n))
  +s\e^{\sup_{\rr\le n+1} V} \mu (1_{\{\rho> n+k+2\}}\GG(fg_n))+\bb(s)\e^{J_{n,k}} \mu_V(|f|)^2 \\
&\le 2s \e^{K_{n,k}}\EE_V(f) +16\ll s\e^{K_{n,k}}\mu_V(f^2) + 16
s\e^{Z_n} \eta_{n,k} \|f\|_\infty^2
 +\bb(s) \e^{J_{n,k}} \mu_V(|f|)^2.\end{split}\end{equation*}
\end{proof}

Now, we are in a position to prove Theorem \ref{finite-sp}.

\beg{proof}[Proof of Theorem $\ref{finite-sp}$] It suffices to prove
for $f\in\A$.

(1) Since $q(x,y)=0$ for $d(x,y)> k_0$, we have
$\gg_{n,k}=\eta_{n,k}=0$ for all $k\ge k_0$ So, by Lemmas \ref{L3.1}
and \ref{L3.2}, \beq\label{DF}\beg{split} \mu_V(f^2)\le
&2\big(6\vv_{n,k} +  s\e^{K_{n,k}}\big)\EE_V(f) + 16\ll\big(8
\vv_{n,k}
 + s\e^{K_{n,k}}\big)\mu_V(f^2)\\
& +\bb(s)\e^{J_{n,k}}\mu_V(|f|)^2,\ \ \ s>0, n\ge 1, k\ge
k_0.\end{split}\end{equation}  If $\inf_{n\ge1,k\ge k_0}\vv_{n,k}=
0$, then for any $r>0$ and any $r'\in (0,r]$ there exist $s>0$,
$n\ge1$ and $k\ge k_0$ such that
$$8\vv_{n,k} +s \e^{K_{n,k}}\le \ff {r'}{2+16\ll r'}.$$ Combining this with (\ref{DF}) we obtain
$$\mu_V(f^2)\le r'\EE_V(f) + (1+8\ll r')\bb(s) \e^{J_{n,k}}\mu_V(|f|)^2\le r\EE_V(f) + (1+8\ll r')\bb(s) \e^{J_{n,k}}\mu_V(|f|)^2.$$ This implies
the super Poincar\'e inequality for the desired $\bb_V$.

(2) According to the expression of the Dirichlet form $(\E_V,
\D(\E_V))$, it is easy to check that $(\E_V, \D(\E_V))$ is
irreducible, i.e. $\E_V(f,f)=0$ implies that $f$ is a constant function.
Indeed, for every $f \in \D(\E_V)$,
\begin{equation*}
\E_V(f,f)=\int_{E\times E}\big(f(x)-f(y)\big)^2q(x,y)\,\mu(\d
y)\,\mu_V(\d x).
\end{equation*}
Hence, if $\E_V(f,f)=0$ for some $f \in \D(\E_V)$, then
$\big(f(x)-f(y)\big)^2=0$ for $\mu_V \times \mu_V$-a.e.\ $(x,y)\in
E\times E $, which implies that $f$ equals to a constant
$\mu_V$-a.e.

According to
   \cite[Corollary 1.2]{W13} (see also \cite[Theorem 1]{Mi}), for the symmetric, conservative, irreducible
 Dirichlet form $(\E_V, \D(\E_V))$,  the desired Poincar\'{e} inequality \eqref{th-sp-2-1} is equivalent to the following defective Poincar\'{e} inequality
\beq\label{th-sp-2}\mu_V(f^2)\le C_1\EE_V(f) +C_2\mu_V(|f|)^2,\ \
f\in\D(\EE_V)\end{equation}  for some constants $C_1,C_2>0$.
Therefore, below we only need to verify \eqref{th-sp-2}.

If \beq \label{AAA} \inf_{n\ge1,k\ge k_0} \vv_{n,k}<\ff 1
{128\ll},\end{equation} then there exist $n\ge 1$, $k\ge k_0$ and
$s>0$ such that
$$16\ll\big(8 \vv_{n,k}
 + s\e^{K_{n,k}}\big)<1.$$ Therefore, the defective Poincar\'e inequality follows from (\ref{DF}).

  To prove (\ref{th-sp-2}) without condition  (\ref{AAA}),    we   follow the approach in the proof of \cite[Theorem 3.1(2)]{BLW} by making bounded perturbations of $V$.  For any $N\ge 1$,
write
$$ V= V_N+ (V\land N)\lor (-N),\ \
V_N:=(V-N)1_{\{V \ge N\}}+(V+N)1_{\{V \le -N\}}.$$ Since $(V\land
N)\lor (-N)$ is bounded and the defective Poincar\'e inequality is
stable under bounded perturbations of $V$, we only need to prove
that when $V$ is unbounded we may find $N\ge 1$ such that   the
defective Poincar\'e inequality holds for $V_N$ in place of $V$.
Note that for any $N\ge 1$, $\mu (e^{V_N})<\infty$. To prove
\eqref{th-sp-2}, without loss of generality, we do not require that
$\mu(e^{V_N})=1$ in the proof below; otherwise, one can easily
normalize it by choosing possibly different positive (but still
finite) constants $C_1$ and $C_2$ in (\ref{th-sp-2}). It is easy to
check that for $n$ large enough
\begin{equation*}
\sup_{\rho \le n}V_N\le
\begin{cases}
& \big(\sup_{\rho \le n}V\big)-N,\ \ \ \ \text{if}\ \ \sup_{E}V=+\infty,\\
& \sup_{\rho \le n}V,\ \ \ \ \ \ \ \ \ \ \quad\ \text{if}\ \
\sup_{E}V<+\infty,
\end{cases}
\end{equation*} and
\begin{equation*}
\inf_{\rho \le n}V_N\ge
\begin{cases}
& \big(\inf_{\rho \le n}V\big)+N,\ \ \ \ \text{if}\ \ \inf_{E}V=-\infty,\\
& \inf_{\rho \le n}V,\ \ \ \ \ \ \ \ \ \ \ \quad \text{if}\ \
\inf_{E}V>-\infty.
\end{cases}
\end{equation*}Thus, for large $n$ we have $K_{n,k}(V_N)\le K_{n,k}(V)-N$, so that
\begin{equation}\label{proof-e-12}
\begin{split}
&\inf_{n\ge1,k\ge k_0}\vv_{n,k}(V_N) \le \e^{-N}\inf_{n\ge 1,k\ge
k_0}\vv_{n,k}(V).
\end{split}
\end{equation}
Since $\inf_{n\ge 1,k\ge k_0}\vv_{n,k}(V)<\infty,$ we see that
(\ref{AAA}) holds for $V_N$ in place of $V$ when $N$ is  large
enough. Therefore, the defective Poincar\'e inequality holds for
$V_N$ in place of $V$ as observed above.
\end{proof}

\subsection{Jumps of infinite range}

When the jump is of infinite range, we will need additional
notations and assumptions to control the uniform norm appearing in
the perturbed functional inequalities (see Lemmas \ref{L3.1} and
\ref{L3.2}), which is an essentially different feature from
the diffusion setting.  For any $n, k\ge 1$ and $\delta>1$, let
\beg{equation*}\beg{split}
 &Z_n(V)= \sup_{\rho(x)\le n+1}V(x),\\
  &\zeta_{n}(V)=\sup_{m\ge n} \Big\{\bb^{-1}\big(1/[2\mu(\rho>m-1)]\big) \e^{Z_{m+1}(V)}\Big\},\\
 &t_{i,n,k}(\dd,V):= \bb^{-1}\Big(\ff 1 4 \dd^i\e^{-J_{n,k}(V)}\Big).\end{split}\end{equation*}
 We assume
 \paragraph{(A)} \ There exist    $\dd>1$ and   sequences $\{(n_i,k_i)\}_{i\ge 1}\subset  \mathbb N^2$   such that $n_i\uparrow\infty$ and\beg{enumerate}
 \item[(A1)] $\lim_{i\to\infty}(\vv_{n_i,k_i}(V)+t_{i,n_i,k_i}(\dd,V)\e^{K_{n_i,k_i}(V)})=0;$
 \item[(A2)] $\sum_{i=1}^\infty\big( \zeta_{n_i}(V)\gg_{n_i,k_i} +t_{i,n_i,k_i}(\dd,V)\e^{Z_{n_i}(V)}\eta_{n_i,k_i}\big)\delta^i <\infty.$\end{enumerate}

 \
Assumptions (A1) and (A2) refer to controls of the perturbations in
terms of  the finite range jumps part and the infinite range jumps
part  respectively. Similarly to the condition in Theorem
\ref{finite-sp}(1), assumption (A1) means that the growth of $V$ is
dominated according to behaviors of $\mu$ and $\bb$  in the  super
Poincar\'{e} inequality $(\ref{sp})$. On the other hand, assumption
(A2) balances the growth condition on $V$ and the intensity of the
  jump kernel $q$ of infinite range,  which is automatically
satisfied if the range of jump is finite   (i.e.\ there exists
$k_0\ge 1 $ such that $q(x,y)=0$ for  $d(x,y)> k_0$) since in the
case $\gg_{n,k}=\eta_{n,k}=0$ for $k\ge k_0$.

 Let $I_\dd$ denote the set of all sequences $\{(n_i,k_i)\}\subset \mathbb N^2$ such that $n_i\uparrow \infty$ and (A1)-(A2) hold. Moreover, for any $r>0$ and $\{n_i,k_i\}\in I_\dd$,
 let $D(r,\{(n_i,k_i)\})$ be the set of $j\in\mathbb N$ such that
 \beg{equation}\label{A-01}\sup_{i\ge j}\big(8\vv_{n_i,k_i}(V)+   t_{i,n_i,k_i}(\dd,V)
 \e^{K_{n_i,k_i}(V)}\big)\le \ff{1}{64\ll}\land \ff{c(\dd)r}{16},\end{equation} where $ c(\dd):= \big(\ff{\ss\dd-1}\dd\big)^2,$ and such that
 \beq\label{A-02} \sum_{i=j}^\infty \Big(6\zeta_{n_i}(V)\gg_{n_i,k_i} +t_{i,n_i,k_i}(\dd,V)\e^{Z_{n_i}(V)}\eta_{n_i,k_i}\Big)\dd^{i+2}\le \ff{1}{256}.\end{equation}
 By {\bf (A)}, we see that for any $r>0$ and $\{(n_i,k_i)\}\in I_\dd$, the set $D(r,\{(n_i,k_i)\})$ is non-empty.

\beg{thm}\label{th-sp} Assume that $(\ref{sp})$ holds.
\begin{itemize}
\item[$(1)$] If {\bf
(A)} is satisfied, then the super Poincar\'e inequality
\eqref{th-sp-1}
 holds with
 $$\bb_V(r):= \inf\bigg\{2\dd^j:\ \{(n_i,k_i)\}\in I_\dd,\ j\in D(r,\{(n_i,k_i)\})\bigg\} <\infty,\ \ r>0.$$

\item[$(2)$] If $(A2)$ is satisfied and $(A1)$ is replaced by the following weaker assumption
\item[$(A1')$]$\limsup_{i\to\infty}(\vv_{n_i,k_i}(V)+t_{i,n_i,k_i}(\dd,V))<\infty,$\newline
then the  Poincar\'e inequality \eqref{th-sp-2-1} holds for some
$C>0.$
\end{itemize}
 \end{thm}

 To prove  Theorem $\ref{th-sp}$, we need to get rid of
the uniform norm included in Lemmas \ref{L3.1} and \ref{L3.2}, for which  we
adopt a cut-off
 argument as in the proof of \cite[Theorem 3.2]{W00a} or \cite[Theorem 3.3.3]{Wbook}. More precisely,  for $\dd>1$ in assumption {\bf (A)} and a  non-negative function $f$, let
 \begin{equation}\label{pro-function}f_{\dd,i}  =(f-\dd^{\ff i 2})^+\land (\dd^{\ff{i+1}2} -\dd^{\ff i 2}),\ \ i\ge 0.\end{equation}
 According to \cite[Lemma 3.3.2]{Wbook}, for $f\in \D(\EE_V)$ and $i,j\ge 0$, we have $f_{\dd,i}$, $(f-\dd^{\ff j 2})^+$ and $f\land \dd^{\ff j 2}\in\D(\EE_V)$. Moreover,
 \beq\label{dri-01} \sum_{i=j}^\infty \EE_V(f_{\dd,i})\le \EE_V((f-\dd^{\ff j 2})^+),\ \ \EE_V((f-\dd^{\ff j 2})^+)+\EE_V(f\land \dd^{\ff j 2})\le \EE_V(f).\end{equation}
We also have the following lemma.

 \beg{lem}\label{L3.3} For any non-negative function $f$ and $k\in \mathbb Z_+:=\{0,1,2,\cdots\}$,
 \beq\label{dri-02} \sum_{i=k}^\infty f_{\dd,i}^2\ge c(\dd) \big({(f-\dd^{\ff k 2})^+}\big)^2.\end{equation} \end{lem}

 \beg{proof} We shall simply use $f$ to denote its value at a fixed point.
 If $ f \le 1$, then both sides in (\ref{dri-02}) are equal to zero. Assume that $f\in (\dd^{\ff l 2}, \dd^{\ff{l+1}2}]$ for some $l\in \mathbb Z_+$. If $l\le k$ then
 $$\sum_{i=k}^\infty f_{\dd,i}^2= \big({(f-\dd^{\ff k 2})^+}\big)^2 \ge c(\dd)\big({(f-\dd^{\ff k 2})^+}\big)^2. $$ Next, if $l>k$ then
 $$\sum_{i=k}^\infty f_{\dd,i}^2\ge (\dd^{\ff l 2}-\dd^{\ff{l-1}2})^2= c(\dd)\dd^{l+1} \ge c(\dd)f^2.$$
 In conclusion, (\ref{dri-02}) holds. \end{proof}

\beg{proof}[Proof of Theorem $\ref{th-sp}$] Since $\EE_V(|f|,|f|)\le
\EE_V(f,f)$ for every
 $f \in \A$, without loss of generality, we may and do assume that $f\in\A$ with $f\ge 0$ and $\mu_V(f^2)=1$.

(1) By Lemmas \ref{L3.1} and \ref{L3.2}, we have
  \beq\label{proof-e-09}\beg{split}\mu_V(f^2)\le &
2\big(6\vv_{n,k} +  s\e^{K_{n,k}}\big)\EE_V(f) + 16\ll\big(8 \vv_{n,k}+ s\e^{K_{n,k}}\big)\mu_V(f^2)\\
&+\bb(s)\e^{J_{n,k}}\mu_V(|f|)^2+16\big(6\zeta_{n}\gg_{n,k}+s\e^{Z_{n}}\eta_{n,k}
\big) \|f\|_\infty^2,\ \ s,k\ge 1,n \ge 1.\end{split}\end{equation}

Next,  let $r>0$, $\{(n_i,k_i)\}\in I_\dd$, $j\in
D(r,\{(n_i,k_i)\})$ be fixed, and let $f_{\dd,i}$ be   defined by
\eqref{pro-function}. Since $\|f_{\dd,i}\|^2_\infty\le
c(\delta)\dd^{i+2}$ and by the Cauchy-Schwarz inequality
$$\mu_V(f_{\dd,i})^2=\mu_V(f_{\dd,i}1_{\{f\ge \delta^{i/2}\}})^2
\le \mu_V(f_{\dd,i}^2)\mu_V(f^2\ge \dd^i)\le
\mu_V(f_{\dd,i}^2)\dd^{-i},$$ it follows from (\ref{A-01}) and
(\ref{proof-e-09}) with $n=n_i$ and $s= t_{i,n_i,k_i}$   that
$$\mu_V(f_{\dd,i}^2)\le \frac{c(\dd)r}{8} \EE_V(f_{\dd,i})
+\ff 1 2 \mu_V(f_{\dd,i}^2) + 16c(\delta)(6
\zeta_{n_i}\gg_{n_i,k_i}+t_{i,n_i,k_i}\e^{Z_{n_i}} \eta_{n_i,k_i})
\dd^{i+2},\ \ i\ge j.$$ That is,
$$\mu_V(f_{\dd,i}^2)\le \frac{c(\dd)r}{4} \EE_V(f_{\dd,i})
+ 32c(\delta)(6 \zeta_{n_i,k_i}\gg_{n_i,k_i}+t_{i,n_i,k_i}
\e^{Z_{n_i}}\eta_{n_i,k_i})\dd^{i+2},\ \ i\ge j.$$ Taking summation
over $i\ge j$ and using (\ref{dri-01}), (\ref{dri-02}) and
(\ref{A-02}), we obtain \beq\label{proof-e-10}
\mu_V\big({(f-\dd^{\ff j 2})^+}^2\big)\le \ff r 4 \EE_V((f-\dd^{\ff
j 2})^+) +\ff 1 8.\end{equation} On the other hand, noting that
$c(\delta)\in(0,1)$ and $\delta>1$, applying (\ref{proof-e-09}) with
$n=n_j$ and $s= t_{j,n_j,k_j}$ to $f \wedge \dd^{\frac{j}{2}}$, and
combining with (\ref{A-01}), (\ref{A-02}), we obtain
\beg{equation*}\beg{split}
\mu_V(f^2\land\dd^j)&\le\frac{c(\delta)r}{8} \EE_V(f\land  \dd^{\ff
j 2}) +\ff 1 4 \mu_V(f^2\wedge\dd^{ j})+ \ff {\dd^j}4 \mu_V(|f|)^2 \\
&\quad+
16( 6 \vv_{n_j,k_j}\gg_{n_j,k_j}+t_{j,n_j,k_j}\e^{Z_{n_j}}\eta_{n_j,k_j})\dd^{j }\\
&\le \ff r {8} \EE_V(f\land  \dd^{\ff j 2}) +\ff 1 2
\mu_V(f^2\wedge\dd^{ j})+ \ff {\dd^j}4 \mu_V(|f|)^2 + \frac{1}{16}.
 \end{split}\end{equation*} Thus,
$$\mu_V(f^2\land\dd^j)\le \ff r 4 \EE_V(f\land  \dd^{\ff j 2}) +\ff{\dd^j} 2 \mu_V(|f|)^2 +\ff 1 8.$$
Combining this with (\ref{proof-e-10}) and using (\ref{dri-01}), we
arrive at \beg{equation*}\beg{split} 1&= \mu_V(f^2)\le
\mu_V\big\{((f-\dd^{\ff j 2})^++(f\wedge
\dd^{\ff j 2}))^2\big\}\\
&\le 2 \mu_V\Big(\big({(f-\dd^{\ff j 2})^+}\big)^2\Big)+ 2  \mu_V(f^2\land\dd^j)\\
&\le \ff r 2\Big(\EE_V((f-\dd^{\ff j 2})^+)+\EE_V(f\land
\dd^{\ff j2})\Big)+ \ff 1 2 +\dd^j \mu_V(|f|)^2\\
&\le \ff r 2 \EE_V(f) +\ff 1 2 +\dd^j
\mu_V(|f|)^2.\end{split}\end{equation*}  Therefore, \beq\label{DDE}
\mu_V(f^2)\le r\EE_V(f) + 2 \dd^j\mu_V(|f|)^2\end{equation} holds
for all $\{(n_i,k_i)\}\in I_\dd$ and $j\in D(r,\{(n_i,k_i)\})$. This
proves (\ref{th-sp-1}) for the desired $\bb_V$.

(2) Similar to the proof of Theorem \ref{finite-sp}(2), we only need
to prove the defective Poincar\'e inequality \eqref{th-sp-2}.  If
\beq\label{DDF} \limsup_{i\to\infty} (8\vv_{n_i,k_i}
+t_{i,n_i,k_i}e^{K_{n_i,k_i}})\le \ff 1 {64\ll},\end{equation}then
there exist a constant $r>0$ and $j\ge 1$
  such that
(\ref{A-01}) and (\ref{A-02}) hold, i.e.\ $j\in D(r,\{(n_i,k_i)\})$.
So, the arguments in (1) ensure    (\ref{DDE}),  and so the
defective Poincar\'e inequality for $C_1=r$ and $C_2= 2\dd^j.$ Let
$V_N$ be in the proof of Theorem \ref{finite-sp}(2). It follows from
the proof of Theorem \ref{finite-sp}(2) that for $n$ large enough
$K_{n,k}(V_N)\le K_{n,k}(V)-N$. Since $J_{n,k}(V_N)\le J_{n,k}(V)$
and $\beta$ is decreasing, we have
\begin{equation*}
\begin{split}
&\limsup_{i \to\infty}t_{i,n_i,k_i}(\dd,V_N)e^{K_{n_i,k_i}(V_N)} \le
\e^{-N}\limsup_{i\to\infty}t_{i,n_i,k_i}(\dd,V)e^{K_{n_i,k_i}(V)},
\end{split}
\end{equation*} which combined with \eqref{proof-e-12} yields that
$$\limsup_{i \to\infty}\Big(\vv_{n_i,k_i}(V_N)+t_{i,n_i,k_i}(\dd,V_N)e^{K_{n_i,k_i}(V_N)} \Big)\le
\e^{-N}\limsup_{i\to\infty}\Big(\vv_{n_i,k_i}(V)+t_{i,n_i,k_i}(\dd,V)e^{K_{n_i,k_i}(V)}\Big).$$
Then the remainder of the proof is similar to that of Theorem
\ref{finite-sp}(2) by using (\ref{DDF}) instead of (\ref{AAA}).
\end{proof}

\subsection{Examples}
The following example shows that Theorem \ref{th-sp} is sharp in
some specific situations.

\begin{exa}\label{ex-sp-1}\rm
Let $E=\R^m$ with $d(x,y)=|x-y|$, and let   $$q(x, y)=\frac{(1+|y|)^{m+\aa}\log(1+|y|)}{|x-y|^{m+\alpha}},\ \ \mu(\d
x)=\frac{c_{m,\aa}\d x}{(1+|x|)^{m+\alpha}\log(1+|x|)},$$   where
$ \alpha\in (0,2)$ and $c_{m,\aa}>0$ is the normalizing constant such that $\mu$ is a probability measure. It is easy
to see that all the assumptions in the introduction for
$(\EE,\scr{D}(\EE))$ are satisfied. We consider $V$ satisfying
\begin{equation}\label{ex-sp-1-01}
-s\varepsilon\log\log(e+|x|)-K \le  V(x)\le (1-s)\varepsilon\log\log(e+|x|)+K,\ \ x\in \R^m
\end{equation}
for some constants $\varepsilon\in(0,1]$, $s \in [0,1]$ and $K \in \R$ such that $\mu(\e^V)=1.$
\begin{itemize}
\item[(1)]  If $\vv<1$ then  (\ref{th-sp-1}) holds
with
\begin{equation}\label{ex-sp-1-02}
\beta_V(r)=\exp\big(C_1(1+{r}^{-1/(1-\vv)})\big),
\end{equation}
for some constant $C_1>0$.
\item[(2)] $\beta_V$ in (1) can not be replaced by any essentially smaller functions, i.e.
 when $V(x)=\vv\log\log(\e+|x|)+K_0$ for some constant $K_0\in \R$ such that $\mu(\e^V)=1$, the estimate
\eqref{ex-sp-1-02} is sharp in the sense that the super Poincar\'e
inequality (\ref{th-sp-1}) does not hold if
$$\lim_{r\to 0}r^{1/(1-\vv)}\log \bb_V(r)=0.$$
\item[(3)] In (1) the constant $K$ can not be replaced by any unbounded positive function, i.e. for
any increasing function $\phi: [0,\infty)\to [0,\infty)$ with $\phi(r)\uparrow\infty$ as $r\uparrow\infty$, there exists $V$ such that \begin{equation}\label{ABC1}
  \vv\log\log (\e+|x|)\le V(x)\le  \varepsilon\log\log(\e+|x|)+\phi(|x|)
\end{equation} with $\mu(\e^V)<\infty,$ but     the super Poincar\'e
inequality (\ref{th-sp-1}) with   $\beta_{V}$
given by (\ref{ex-sp-1-02}) does not hold for $\mu_V(\d x):= \ff{\e^{V(x)}\d x}{\mu(\e^V)}.$
\item[(4)] If $\vv=1$, then
$(\EE_V,\scr{D}(\EE_V))$ satisfies the  Poincar\'e
inequality
\begin{equation}\label{ex2.2}
\mu_V(f^2)\le C\EE_V(f,f)+\mu_V(f)^2,\ \ \ f \in \scr{D}(\EE_V)
\end{equation} with some constant $C>0$.

\end{itemize}
\end{exa}

\begin{proof} As (2) is included in \cite[Corollary 1.3]{WW}, we only prove (1), (3) and (4).

(a) According to \cite[Corollary 1.3(3)]{WW}, we know the
logarithmic Sobolev inequality holds for $(\EE,\scr{D}(\EE))$, i.e.\
the super Poincar\'{e} inequality \eqref{sp} holds for
$(\EE,\scr{D}(\EE))$ with the rate function
$\beta(r)=\exp\big(c_1(1+{r}^{-1})\big)$ for some constant $c_1>0$.

Next, by   (\ref{ex-sp-1-01}), there exists a constant $c_2>0$
such that for $n$ large enough
\begin{equation}\label{e3.10}
\begin{split}
& K_{n,n}(V)\le \vv\log\log n+c_2,\ \ J_{n,n}(V)\le (1+s)\vv \log\log n+c_2,\ \\
& Z_{n}(V)\le (1-s)\vv \log\log n+c_2,\ \ \zeta_n(V)\le \frac{c_2}{\log^{1-(1-s)\vv}n},\ \ \vv_{n,n}(V)\le
\frac{c_2}{\log^{1-\vv} n},\\
&\gg_{n,n}\le \frac{c_2}{n^{\alpha}},\ \ \eta_{n,n}\le \frac{c_2}{n^{\alpha}\log(1+n)}.
\end{split}
\end{equation}
For any $\delta>1$, taking
$\{n_i\}=\{\delta^i\}_{i=1}^{\infty}$, and $k_i=n_i$, we
get that $t_{i,n_i,k_i}(\dd,V)\le {c_3}{i}^{-1}$ for $i$ large
enough and some constant $c_3>0$. So, assumption {\bf (A)} is
fulfilled.

Moreover, by (\ref{A-01}) and (\ref{A-02}), it is easy to check that
for $r>0$ small enough we have
\begin{equation*}
\Big\{j \ge1:  j \ge c_4r^{-\frac{1}{1-\vv}}\Big\}\subseteq D(r,\{(n_i,k_i)\}).
\end{equation*}
Then, the   assertion in (1)  follows from Theorem \ref{th-sp}(1) and \eqref{e3.10}.

(b) Let $ \psi(r)= (1+\log(1+r))\land \e^{\phi(r)}$ for $r\ge 0.$
Then $1\le\psi\le \e^\phi, \psi(r)\uparrow\infty$ as
$r\uparrow\infty$ and $\psi(r)\le 2\log r$ for large $r$. Let
$$V(x)= \vv\log\log (\e+|x|) +\log\psi(|x|).$$ Then $V$ satisfies (\ref{ABC1}) and $\mu(\e^V)<\infty.$ Up to a normalization constant we may simply assume that $\mu_V(\d x)= \e^{V(x)}\d x.$

Now, suppose that  (\ref{th-sp-1}) holds  with $\bb_V$ given by
(\ref{ex-sp-1-02}).
For any $n\ge1$, let $f_n\in C^\infty(\R^d)$ satisfy that
\begin{equation*} f_{n}(x)
\begin{cases}
& =0,\ \ \ \ \ \ \ \quad \,\,|x|\le n,\\
& \in [0,1],\ \ \ \ \quad n\le |x|\le 2n,\\
& =1,\ \ \ \ \ \ \ \ \quad |x|\ge 2n,
\end{cases}
\end{equation*} and $|\nabla f_n|\le \frac{2}{n}.$ Then, there exists a constant $c_5>0$ independent of $n$ such that
\begin{equation}\label{AB0}
\begin{split}
\Gamma(f_n,f_n)(x)&=c_{m,\aa}\int\frac{(f_n(y)-f_n(x))^2}{|y-x|^{m+\aa}}\,\d y\\
&\le \frac{4c_{m,\aa}}{n^2}\int_{\{|y-x|\le n\}}\frac{\d y}{|y-x|^{m+\alpha-2}}+c_{m,\aa}\int_{\{|y-x|> n\}}\frac{\d y}{|y-x|^{m+\alpha}}\\
&\le \frac{c_5}{n^\alpha},\quad n\ge1.
\end{split}
\end{equation}
According to the definition of $f_n$ and the increasing property of $\psi$, there exists $c_6>0$ such that
\begin{equation}\label{proof-e-0100}\mu_V(f_n^2)\ge \frac{c_6  \psi(n)}{n^{\aa}\log^{1-\vv}n},\ \ n\ge 2.\end{equation}
On the other hand, since $\psi(r)\le 2\log r$ for large $r$, we have
$$\mu_V(|\cdot|>n)\le c_7 \int_n^\infty \ff{\psi(r)\d r}{r^{1+\aa}\log^{1-\vv}r }\le
\ff{c_8\log^{\vv} n}{n^\aa},\ \ n\ge 2.$$ Then
$$\mu_V(f_n)^2\le \mu_V(f_n^2) \mu_V(|\cdot|>n) \le \ff{c_8\log^{\vv}n}{n^\aa}\mu_V(f_n^2),\ \ n\ge 2.$$
Combining this with (\ref{th-sp-1}), (\ref{AB0}) and
(\ref{proof-e-0100}), we obtain
\beq\label{AD}1-\ff{c_8\log^{\vv}n}{n^\aa}\bb_V(r)\le\ff{c_5r}{n^\aa\mu_V(f_n^2)}
\le \ff{c_9r\log^{1-\vv}n}{\psi(n)},\ \ r>0, n\ge 2\end{equation}
for some constant $c_9>0.$ Since $\psi(n)\le 2\log n$ for large $n$,
it is easy to see that
$$r_n:= \bb_V^{-1}\Big(\ff{n^\aa}{2c_8\log^{\vv}n}\Big)\le \ff{c_{10}}{\log^{1-\vv}n}$$ holds for some constant $c_{10}>0$ and large $n$. Therefore, it follows from (\ref{AD}) that
$$\ff 1 2 \le \lim_{n\to\infty} \ff{c_9r_n\log^{1-\vv}n}{\psi(n)}=0,$$ which is a contradiction.

(c) If $\vv=1$ in (\ref{ex-sp-1-01}), the estimate (\ref{e3.10}) is
still true,
and so the required Poincar\'e inequality follows from Theorem \ref{th-sp}(2). \end{proof}

Finally,  we consider an example for finite range of jumps to
illustrate Theorem \ref{finite-sp}.

\begin{exa}\label{ex2}\rm
Let $E=\R^m$  with $d(x,y)=|x-y|$ and let
$$q(x, y)=\frac{\e^{|y|^\kk}}{{|x-y|}^{m+\alpha}}1_{\{{|x-y|}\le 1\}} ,\ \  \mu(\d x)=c_{m,\kappa}\e^{-|x|^{\kappa}}\d x $$ for some
constants $0<\alpha<2$,  $\kappa>1$ and $c_{m,\kk}\ge 1$ such that
$\mu$ is a probability measure.
 It is easy to see that all the assumptions in the introduction for
$(\EE,\scr{D}(\EE))$ is satisfied. Consider $V$ satisfying
\begin{equation}\label{ex2.0}
-C_1(1+|x|^{\theta-1})-K \le  V(x)\le C_1(1+|x|^{\theta-1})+K,\quad
x\in \R^m
\end{equation}
for some constants $\theta \in (1,\kappa]$, $C_1>0$ and $K\in\R$
such that $\mu(\e^V)=1.$
\begin{itemize}
\item[(1)] If $\theta<\kk,$ then   the super Poincar\'e inequality holds
for $(\EE_V, \scr{D}(\EE_V))$ with
$$
\beta_V(r)=C_2\exp\Big[C_2 \log^{\frac{\kappa}{\kappa-1}}(1+r^{-1})
\Big]
$$
for some positive constants $C_2>0.$

\item[(2)] Let $\theta=\kk$. Then if $C_1>0$ is small enough,  then   the Poincar\'e inequality \eqref{ex2.2}  holds for some constants $C>0.$
\end{itemize}
\end{exa}

\begin{proof} According to \cite[Example 1.2(2)]{CW2}, we know the super Poincar\'{e} inequality
\eqref{sp} holds for $(\EE,\scr{D}(\EE))$ with
\beq\label{DG0}\beta(r)=c_1\exp\Big[c_1\log^{\frac{\kappa}{\kappa-1}}(1+
r^{-1})\Big]\end{equation} for some constant $c_1>0$.

(1) By (\ref{ex2.0}) and the fact that $\theta<\kappa$, one can find
some constants $c_2, c_3>0$ such that for $n$ large enough
\begin{equation*}
\begin{split}
& K_{n,n}(V)\le c_2n^{\theta-1},\ \ J_{n,n}(V)\le c_2n^{\theta-1},\
\ \vv_{n,n}(V)\le \e^{-c_3n^{\kappa-1}}.
\end{split}
\end{equation*}
Therefore, for every $r$ small enough, taking
$$n=c_4\log^{\frac{1}{\kappa-1}}(1 +{r}^{-1})$$ and
$$s=\exp\big(-c_5\log(1+{r}^{-1})\big)$$ for some constants
$c_5>>c_4>>1$,  we obtain
\begin{equation*}
8\vv_{n,n}(V)+s\e^{K_{n,n}(V)}\le \frac{r}{2(1+8\ll r)}.
\end{equation*}
The first required assertion for super Poincar\'{e} inequality of
$(\EE_V,\scr{D}(\EE_V))$ follows from Theorem \ref{finite-sp}(1) and
all the estimates above.

(2) Let $\theta=\kk$. By \eqref{ex2.0}, (\ref{DG0}) and the definition of $\mu$,    for $n$ large enough we have
\begin{equation*}
K_{n,n}(V)\le 3C_1n^{\kappa-1},\ \ \bb^{-1}(\mu(|x|>n-1)^{-1})\le
\e^{-c_6n^{\kk-1}},
\end{equation*} for some constant $c_6>0$ depending only on $\bb$ and $\mu$.
So, if $C_1< c_6/3$, then $\inf_{n\ge1}\vv_{n,n}(V)<\infty$.
Therefore,  by Theorem \ref{finite-sp}(2), the  Poincar\'e
inequality holds for $(\EE_V, \scr{D}(\EE_V))$. \end{proof}

To show that in  Example \ref{ex2}(2) it is essential to assume that
$C_1>0$ is small, we present below a counterexample inspired by
\cite[Proposition 5.1]{BLW}.

\beg{prp} In the situation of Example $\ref{ex2}$, let $\theta=\kk, \aa\in (0,1)$ and $m=1.$  Let
$$V(x)= K_0+ L \sum_{n=1}^\infty 1_{[nH, (n+1)H)}(x)(n+1)^{\kk-1} \Big( 2n+1-\ff{2x}H\Big),$$
where $H>4, L>\ff{\kk H^\kk}{H-2}$ and $K_0\in \R$ are constants
such that $\mu(\e^V)=1.$ Then   $(\ref{ex2.0})$ holds for some
constant $C_1>0$ and $K\in\R$; however, for any $C>0$,  the
Poincar\'e inequality $(\ref{ex2.2})$ does not hold.\end{prp}

\beg{proof}   It suffices to disprove the Poincar\'e inequality, for which we are going to construct a sequences of functions $\{f_n\}\subset\A$ such that
\begin{equation*}\label{DPP} \lim_{n\to\infty} \ff{\EE_V(f_n)}{{\rm Var}_{\mu_V} (f_n)}=0.\end{equation*}
For $n\ge 1$, let
\begin{equation*}
f_{n}(x)=\begin{cases}
& \frac{\int_{Hn+1}^{x}\exp(y^{\kappa}-V(y))\d y}
{\int_{Hn+1}^{H(n+1)-1}\exp(y^{\kappa}-V(y))\d y},\ \ \ x \in (Hn+1,H(n+1)-1) ,\\
& \qquad\quad1,\ \ \ \ \ \ \qquad\qquad \,\,\,\,x \in [H(n+1)-1,+\infty),\\
& \qquad\quad0,\ \ \ \ \ \ \ \ \qquad \qquad\text{otherwise},
\end{cases}
\end{equation*}
which is a bounded Lipschitz continuous function on $\R$, so that $f_n\in\A$. In the following calculations $C$ stands for a constant which   varies from line to line but is independent of $n$ (may depend on $H$ or $L$). We simply denote $K_n=L(n+1)^{\kk-1}$ for $n\ge 1$, so that
$$V(x)=K_0+ \sum_{n=1}^\infty 1_{[nH, (n+1)H)}(x)K_n\Big(2n+1-\ff{2x}H\Big).$$

(a) Estimate on ${\rm Var}_{\mu_V}(f_n).$ For $n$  large enough, since
$z\mapsto\big(z^{\kappa}-K_{n+1}(2n+3-\ff{2z}H)\big)$ is increasing, we have
\begin{equation*}
\begin{split}
\mu_{V}(f_{n}^2)  &\ge \mu_{V}\big(H(n+1)\le x \le H(n+2)\big)\\
&\ge C\int_{H(n+1)}^{H(n+2)}\frac{\exp\big(K_{n+1}(2n+3-\ff{2z}H)-z^{\kappa} \big)}{\kappa z^{\kappa-1}+\frac{2K_{n+1}}{H}} \d\bigg(z^{\kappa}-K_{n+1}(2n+3-\ff{2z}H)\bigg)\\
&\ge \frac{C\big\{\exp\big(-(H(n+1))^{\kappa}+K_{n+1}\big)-\exp\big(-(H(n+2))^{\kappa}-K_{n+1}\big)\big\}}{ (n+2)^{\kappa-1}}\\
&\ge \frac{C\exp\big(-(H(n+1))^{\kappa}+K_{n+1}\big)}{ (n+2)^{\kappa-1}}.
\end{split}
\end{equation*}  Noting that for $n$ large enough,
\begin{equation*}
\mu_{V}(f_{n})^2=\mu_{V}(f_{n} 1_{\{x \ge n\}})^2\le \mu_{V}(f_{n}^2)\mu_{V}(x \ge n)
\le \frac{ \mu_{V}(f_{n}^2)}{2},
\end{equation*} we arrive at
\beq\label{PRO-01} {\rm Var}_{\mu_V}(f_n)\ge
\frac{C\exp\big(-(H(n+1))^{\kappa}+K_{n+1}\big)}{(n+2)^{\kappa-1}}.\end{equation}

(b) Estimate on $\EE_V(f_n)$. Let $$g_{n}(x) =\bigg(\int_{Hn+1}^{H(n+1)-1}\exp\big(y^{\kappa}-V(y)\big)\d y\bigg)f_{n}(x).$$
Noting that
$g_{n}'(x)\neq 0$ only when $x \in (Hn+1,H(n+1)-1)$,
we have
\begin{equation*}
\begin{split}
\EE_{V}(g_{n})& \le C\int_{Hn}^{H(n+1)}\bigg(\int_{\{|x-y|\le 1\}}
\frac{|\int_{x}^{y} g_{n}'(z)\d z|^2}{{|x-y|}^{1+\alpha}}\d y \bigg)\exp\big(-x^{\kappa}+V(x)\big) \d x\\
&=C\int_{Hn}^{H(n+1)}\bigg(\int_{\{{|x-y|}\le 1, x \le y\}}
\frac{|\int_{x}^{y} g_{n}'(z)\d z|^2}{{|x-y|}^{1+\alpha}}\d y\bigg)\exp\big(-x^{\kappa}+V(x)\big) \d x\\
&\quad+C\int_{Hn}^{H(n+1)}\bigg(\int_{\{{|x-y|}\le 1, x > y\}}
\frac{|\int_{x}^{y} g_{n}'(z)\d z|^2}{{|x-y|}^{1+\alpha}}\d y \bigg)\exp\big(-x^{\kappa}+V(x)\big) \d x\\
&=:C(I_1+I_2).
\end{split}
\end{equation*}
Since
$g_{n}'(x)=\exp\big(x^{\kappa}-V(x)\big)$ for $x \in (Hn+1,H(n+1)-1)$, and since for large $n$ the function $z\mapsto z^\kk-V(z)$ is increasing on $(Hn, H(n+1))$,
by the Cauchy-Schwarz inequality we have, for $x\in (Hn, H(n+1)-1)$ and $ x\le y\le x+1$,
\beg{equation*}\beg{split} & \frac{|\int_{x}^{y} g_{n}'(z)\d z|^2}{{|x-y|}^{1+\alpha}}  \\
&=  \frac{|\int_{x \vee (Hn+1)}^{y\wedge (H(n+1)-1)} g_{n}'(z)\d z|^2}{{|x-y|}^{1+\alpha}} \\
&\le \frac{|\int_{x \vee (Hn+1)}^{y\wedge (H(n+1)-1)}
\exp\big(z^{\kappa}-V(z)\big)\d z|}{{|x-y|}^{1+\alpha}}\int_{x \vee (Hn+1)}^{y\wedge (H(n+1)-1)} |g_{n}'|^2(z)\e^{-z^{\kappa}+V(z)}\d z \\
&\le \frac{\exp\big((x+1)^{\kappa}-V(x+1)\big)}{{|x-y|}^{\alpha}}
\Big(\int_{Hn+1}^{H(n+1)-1}|g_{n}'|^2(z)\e^{-z^{\kappa}+V(z)}\d
z\Big).
\end{split}\end{equation*} Thus, since $\aa\in (0,1)$,
\begin{equation*}
\begin{split}
I_1 &\le
C\Big(\int_{Hn+1}^{H(n+1)-1}|g_{n}'|^2(z)\exp\big(-z^{\kappa}+V(z)\big)\d
z\Big) \\
&\qquad\times\bigg[\int_{Hn}^{H(n+1)-1}\Big(\int_{\{|x-y|\le 1\}}
\frac{{1}}{{|x-y|}^{\alpha}}\d y\Big)\exp\Big((x+1)^{\kappa}-x^{\kappa}-
V(x+1)+V(x)\Big) \d x\bigg]\\
&\le C  \Big(\int_{Hn+1}^{H(n+1)-1}\exp\big(z^{\kappa}-V(z)\big)\d z\Big)
\exp\Big(\kappa (H(n+1))^{\kappa-1}+\frac{2K_{n}}{H}\Big).
\end{split}
\end{equation*}
 Similarly, we have
\begin{equation*}
I_2 \le C\int_{Hn+1}^{H(n+1)-1}\exp\big(z^{\kappa}-V(z)\big)\d z.
\end{equation*}
So, for $n$ large enough,
\begin{equation*}
\EE_{V}(g_{n}) \le C
\Big(\int_{Hn+1}^{H(n+1)-1}\exp\big(z^{\kappa}-V(z)\big)\d z\Big)
\exp\bigg(\kappa (H(n+1))^{\kappa-1}+\frac{2K_{n}}{H}\bigg).
\end{equation*}
Moreover, for large $n$,
\begin{equation*}
\begin{split}
& \int_{Hn+1}^{H(n+1)-1}\exp\big(z^{\kappa}-V(z)\big)\,\d z\\
&=\int_{Hn+1}^{H(n+1)-1}\frac{\exp\big(z^{\kappa}-\frac{2K_{n}(Hn-z)}{H}-K_{n}\big)}{\kappa z^{\kappa-1}+\frac{2K_{n}}{H}}
 \d\Big(z^{\kappa}-\frac{2K_{n}(Hn-z)}{H}-K_{n}\Big)\\
&\ge \frac{C\exp\big((H(n+1)-1)^{\kappa}+\frac{(H-2)K_{n}}{H}\big)-\exp\big((Hn+1)^{\kappa}-\frac{(H-2)K_{n}}{H}\big)}
{ (n+1)^{\kappa-1}}\\
&\ge \frac{C\exp\big((H(n+1)-1)^{\kappa}+\frac{(H-2)K_{n}}{H}\big)}{ (n+1)^{\kappa-1}},
\end{split}
\end{equation*}
Therefore, for large $n$,
\beq\label{PRO-02}\beg{split}  \EE_V(f_n)&=\ff{\EE_V(g_n)}{(\int_{Hn+1}^{H(n+1)-1}\exp\big(z^{\kappa}-V(z)\big)\d z)^2}\\
&\le \ff{C \exp\big(\kappa (H(n+1))^{\kappa-1}+\frac{2K_{n}}{H}\big)}{\int_{Hn+1}^{H(n+1)-1}\exp\big(z^{\kappa}-V(z)\big)\d z}\\
&\le C(n+1)^{\kk-1} \exp\Big(\kk(H(n+1))^{\kk-1} -(H(n+1)-1)^\kk -\ff{(H-4)K_n}H\Big). \end{split}\end{equation}

(c) Combining (\ref{PRO-01}) with (\ref{PRO-02}) and noting that
$$(H(n+1))^\kk-(H(n+1)-1)^\kk\le\kk (H(n+1))^{\kk-1},$$ we obtain
\begin{equation*}
\begin{split}
&\lim_{n\to\infty} \frac{\EE_{V}(f_{n},f_{n})}{{\rm Var}_{\mu_{V}}(f_{n})}\\
&\le C\lim_{n\to\infty} (n+2)^{2(\kk-1)} \exp\Big(2\kk(H(n+1))^{\kk-1} -\ff{2(H-2)}H K_n\Big) =0\end{split}\end{equation*}
since $$L>\ff{\kk H^\kk}{H-2}.$$
\end{proof}

\subsection{A variation condition
 }

 It is well known that in the diffusion case the super Poincar\'e inequality is stable under Lipschitz perturbations (see \cite[Proposition 2.6]{W00a}). The aim of this section is to extend this result to the non-local setting using a variation condition on ${\rm supp}\,q:=\{(x,y): q(x,y)>0\}$.

 \beg{thm}\label{th-sup} Assume that $(\ref{sp})$ holds. If there exists a constant
  $\kk_1>0$ such that $\kk_2:= \mu_V(\e^{-2V})=\mu(e^{-V})<\infty$ and
 \beq\label{th-sup-1}   |V(x)-V(y)|\le \kk_1(1\land d(x,y)),\ \ (x,y)\in{\rm supp}\,q.\end{equation} Then  $(\ref{th-sp-1})$ holds for
 $$\bb_V(s):= \inf\Big\{ 16\kk_2\bb(r)^3(4+\ll \kk_1^2 s'):\ s'\in (0,s], 0<r\le \ff{s'\e^{-\kk_1}}{4+\ll \kk_1^2s'  }\Big\},\ \ s>0.$$\end{thm}

 \beg{proof}
To prove (\ref{th-sp-1}) for the desired $\bb_V$, we may and do assume that $V$ is bounded. Indeed, for any $n\ge 1$ let
 $$V_n= (V\land n)\lor (-n) -\log \mu(\e^{(V\land n)\lor(-n)}).$$ Then $\mu(\e^{V_n})=1$, (\ref{th-sup-1}) holds for $V_n$ in place of $V$, and
 $$\lim_{n\to\infty} \mu_{V_n}(\e^{-2V_n})=\kk_2.$$ Thus, applying the assertion to the bounded  $V_n$   and letting $n\to\infty$, we complete the proof.

 Now, let $V$ be bounded and let $f\in\A$ with $\mu_V(|f|)=1$. Take $\tt f=f\e^{\ff V2}.$
 By (\ref{th-sup-1}) we have for every $x,y \in \text{supp}q$,
 \begin{equation*}
 \begin{split}
 &\big(1-\e^{\frac{V(x)-V(y)}{2}}\big)^2 \le \frac{\e^{\kk_1}\kk_1^2 }{4}\big(1\wedge d^2(x,y)\big),
 \end{split}
 \end{equation*}
 hence
 \beg{equation*}\beg{split} \GG(\tt f)(x)&\le 2\int_{E} \Big\{\e^{V(y)}(f(x)-f(y))^2+ f^2(x)(\e^{\ff {V(x)}2}-\e^{\ff {V(y)}2})^2\Big\}q(x,y)\mu(\d y)\\
 &\le 2\e^{V(x)+\kk_1}\GG(f)(x) + \ff{1}2 \kk_1^2\ll \e^{V(x)+\kk_1}f^2(x).
\end{split}\end{equation*} Since $V$ is bounded,  this implies $\tt f\in\A$. Moreover,    combining this with  (\ref{sp}), we obtain for every $r>0$
 \beg{equation}\label{th-sup-2}\beg{split} \mu_V(f^2)&= \mu({\tt f}^2)\le r \mu(\GG(\tt f)) +\bb(r) \mu(|\tt f|)^2\\
 &\le 2 r \e^{\kk_1}\EE_V(f) + r\ff{  \kk_1^2 }2\ll \e^{\kk_1} \mu_V(f^2) +\bb(r) \mu_V( |f|\e^{-\ff V 2})^2.\end{split}\end{equation}
 Since $\mu_V(|f|)=1$, for any $R>0$ we have $\mu_V(|f|>R)\le \ff 1 R$, and hence,
 \beg{equation*}\beg{split} \mu_V(|f|\e^{-\ff V 2})^2 &\le 2 \mu_V(|f|\e^{-\ff V 2}1_{\{|f|\le R\}})^2 + 2 \mu_V(|f|\e^{-\ff V 2}1_{\{|f|\ge R\}})^2\\
 &\le 2 R \mu_V(|f|^{\ff 1 2} \e^{-\ff V 2})^2 +2 \mu_V(f^2) \mu_V(\e^{-V} 1_{\{|f|>R\}})\\
 &\le 2 R \mu_V(|f|)\mu_V(\e^{-V}) + 2\mu_V(f^2) \ss{\mu_V(|f|>R)\mu_V(\e^{-2V})}\\
  &\le 2 R + 2 \mu_V(f^2)\ff{\ss{\kk_2}}{\ss R},\end{split}\end{equation*}where in the second and the third inequalities we have used the Cauchy-Schwarz inequality. Taking $R= 16\kk_2\bb(r)^2$, we get
  $$\bb(r)\mu_V(|f|\e^{-\ff V 2})^2 \le 32\kk_2\bb(r)^3 +\ff 1 2\mu_V(f^2).$$
  Substituting this into \eqref{th-sup-2}, we arrive at
   $$\mu_V(f^2)\le 4r\e^{\kk_1}\EE_V(f)+ r \kk_1^2 \ll \e^{\kk_1}\mu_V(f^2) +64 \kk_2\bb(r)^3.$$
  Therefore,  for any $s>0$ and $s'\in (0,s]$ such that
  $$r\le \ff {s'\e^{-\kk_1}}{4 + \ll \kk_1^2  s'},$$ we have
  $$\mu_V(f^2)\le s' \EE_V(f)+ 16\kk_2\bb(r)^3(4+\ll\kk_1^2 s')\le s \EE_V(f) + 16\kk_2\bb(r)^3(4+\ll\kk_1^2 s').$$
  This implies the desired super Poincar\'e inequality.
 \end{proof}

Let $E=\R^m$ and $d(x,y)=|x-y|$. If the jumps are of finite range, i.e. there is a constant $k\ge 1$ such that $q(x,y)=0$ for $|x-y|>k$, then (\ref{th-sup-1}) holds for any Lipschitz function  $V$. Therefore, the above theorem implies that the super Poincar\'e inequality is stable for  all Lipschitiz perturbations  as it is known in the diffusion case. Note that,  the defective log-Sobolev inequality
$$ \mu(f^2\log f^2)\le C_1 \EE(f)+C_2,\ \ f\in\D(\EE),\mu(f^2)=1$$ holds for some $C_1,C_2>0$ if and only if
 the super Poincar\'e inequality $(\ref{sp})$ holds for $\bb(r)=\e^{c(1+r^{-1})}$ for some $c>0$, see \cite[Corollary 1.1]{W00b} for $\dd=1$. On the other hand, since the symmetric conservative Dirichlet form $(\E_V, \D(\E_V))$ is irreducible, also by \cite[Corollary 1.2]{W13} (see also \cite[Theorem 1]{Mi}), the defective log-Sobolev inequality above implies the (true) log-Sobolev inequality
$$ \mu(f^2\log f^2)\le C \EE(f),\ \ f\in\D(\EE),\mu(f^2)=1$$ for some constant $C>0$. Therefore,
 we conclude from Theorem \ref{th-sup} that the log-Sobolev inequality is stable under perturbations of Lipschitz functions $V.$  See \cite[Example 1.2]{CW2} for examples of $\mu$ and $q$ having finite range of jumps such that the log-Sobolev inequality holds.

\section{Perturbations for the weak Poincar\'{e} inequality}%\label{section3}

Suppose that the  weak Poincar\'{e} inequality
\begin{equation}\label{wp} \mu(f^2)\le \beta(r)
\EE(f,f)+r\|f\|_\infty^2,\quad r>0, f\in \D(\EE),
\mu(f)=0,\end{equation} holds for some decreasing function
$\beta:(0,\infty)\to (0,\infty).$ To derive the weak Poincar\'e
inequality for $\EE_V$ using growth conditions on $V$, for any
$n,k\ge 1$ let
\begin{equation*}
\begin{split}
&\tt K_{n,k} (V)=\sup_{\rho\le n}V -\inf_{\rho \le n+k+1}V,\quad\quad  \tt Z_n(V)=\sup_{\rho\le n}V,\\
&\tt\eta_{n,k} =\int_{\{\rr(x)> n+k+1, \rr(y)\le n+1\}}q(x,y)\mu(\d
y)\mu(\d x),\ \ \tt \gg_k=\int_{\{d(x,y)>k\}}q(x,y)\mu(\d y)\mu(\d
x).
\end{split}
\end{equation*} It is clear that $\tt\eta_{n,k}\le  \tt \gg_k$. By \eqref{intro-02} we have $\tt\eta_{n,k}\downarrow 0$ as $n\uparrow\infty$ or $k\uparrow\infty$.

\beg{thm}\label{th-wp} Assume that the weak Poincar\'{e} inequality
\eqref{wp} holds. If for any $\vv>0$ \beq\label{th-wp-1}\inf_{n,k\ge
1} \e^{\tt Z_{n}(V)}\beta\big(\vv \e^{-\tt
Z_{n}(V)}\big)\big(\tt\eta_{n,k} +\tt
\gg_k+\mu(\rr>n-k)\big)=0,\end{equation} then
 \begin{equation}\label{th-wp-2} \mu_V(f^2)\le
\beta_V(r)\EE_V(f,f)+r\|f\|_\infty^2,\quad r>0,f\in \D(\EE_V),
\mu_V(f)=0\end{equation} holds for \beg{equation*}\beg{split}
\beta_V(r):=\inf\Big\{&2\beta\Big(\ff r 8 \e^{-\tt
Z_{n}(V)}\Big)\e^{\tt K_{n,k}(V)}:\
 6\mu_V(\rr>n)\\
  &+2\e^{\tt Z_{n}(V)}\beta\Big(\ff r 8 \e^{-\tt Z_{n}(V)}\Big)\big(4\tt \eta_{n,k} +\tt \gg_k+4\ll\mu(\rr>n-k)\big)\le \ff r 2\Big\}<\infty,\ \ r>0.
 \end{split}\end{equation*}  \end{thm}

In finite range case, e.g.\ $q(x,y)=0$ for  $d(x,y)\ge1$,
\eqref{th-wp-1} is reduced into $$\inf_{n,k\ge 1} \e^{\tt
Z_{n}(V)}\beta\big(\vv \e^{-\tt Z_{n}(V)}\big)\mu(\rr>n-k)\big)=0,$$
then Theorem \ref{th-wp} is similar to \cite[Proposition 4.1]{BLW}
for local Dirichlet forms.

\beg{proof}[Proof of Theorem $\ref{th-wp}$] It is easy to see that
(\ref{th-wp-1}) implies $\beta_V(r)<\infty$ for $r>0.$ Let $g_n$ be
as in the proof of Lemma \ref{L3.2}. Then for any $f\in\A$ with
$\mu_V(f)=0$, we have
$$\mu_V(fg_{n})^2=\mu_V(f(1-g_{n}))^2\le\|f\|_\infty^2 \mu_V(\rr>n)^2.$$ Moreover,
\beg{equation*}\beg{split} {\rm Var}_{\mu_V}(fg_{n}) &\le \mu_V\big((fg_{n}-\mu(fg_{n}))^2\big)\\
&\le \e^{\sup_{\rr\le n}V}\mu\big(1_{\{\rr\le n\}}(fg_{n}-\mu(fg_{n}))^2\big) +\mu_V\big(1_{\{\rr>n\}}(fg_{n}-\mu(fg_{n}))^2\big)\\
&\le \e^{\sup_{\rr\le n}V}{\rm Var}_{\mu}(fg_{n}) +4\|f\|_\infty^2 \mu_V(\rr>n).\end{split}\end{equation*} Then
\beg{equation}\label{th-wp-3}\beg{split} \mu_V(f^2) &\le {\rm Var}_{\mu_V }(fg_{n})+\mu_V(fg_{n})^2 +\mu_V(f^21_{\{\rr>n\}})\\
&\le \e^{\sup_{\rr\le n} V}{\rm Var}_{\mu}(fg_{n}) +6\|f\|_\infty^2 \mu_V(\rr>n).\end{split}\end{equation}
On the other hand, we have
\beg{equation*}\beg{split} \e^{\sup_{\rr\le n}V} \EE(fg_{n})
  \le &2 \e^{\sup_{\rr\le n}V}\int_{E\times E} g_{ n}^2(y)(f(x)-f(y))^2q(x,y)\mu(\d y)\mu(\d x)\\
  & + 2 \e^{\sup_{\rr\le n}V}\int_{E\times E}f^2(x)(g_{n}(x)-g_{n}(y))^2 q(x,y)\mu(\d y)\mu(\d x)\\
\le &2 \e^{\tt K_{n,k}(V)}\EE_V(f) +2\e^{\tt
Z_n(V)}(4\tt\eta_{n,k}+\tt
\gg_k+4\ll\mu(\rr>n-k))\|f\|_\infty^2,\end{split}\end{equation*}
where in the last inequality we make the first integral on the sets
$\{\rr(x)\le n+k+1\}$ and $\{\rr(x)>n+k+1\}$, and the second
integral on the sets $\{\rr(x)\le n-k\}$ and $\{\rr(x)>n-k\}$, and
also use the facts that \beg{equation*}\beg{split} &\int_{\{\rr(x)>
n+k+1\}}g_n^2(y)(f(x)-f(y))^2
q(x,y)\mu(\d y)\mu(\d x) \le 4\tt\eta_{n,k} \|f\|_\infty^2,\\
&\int_{\{\rr(x)\le n-k\}}f^2(x)(g_n(x)-g_n(y))^2q(x,y)\mu(\d
y)\mu(\d x)
 \le  \tt \gg_k \|f\|_\infty^2,\\
&\int_{\{\rr(x)> n-k\}}f^2(x)(g_n(x)-g_n(y))^2q(x,y)\mu(\d y)\mu(\d
x)\le 4\lambda\mu(\rr>n-k)\|f\|_\infty^2.
\end{split}\end{equation*}Here, in the second inequality we have used the fact that for any $x\in E$ with $\rho(x)\le n-k$, $g_n(x)=1$, and so, by the definition of $g_n$, $|g_n(x)-g_n(y)|\neq0$ only for $y\in E$ with $\rho(y)>n.$
Combining this with (\ref{th-wp-3}) and \eqref{wp}, we arrive for any $s>0$ at
\beg{equation*}\beg{split} \mu_V(f^2) \le &2\beta(s)\e^{\tt K_{n,k}(V)} \EE_V(f)\\
&+ \|f\|_\infty^2\big\{6\mu_V(\rr>n) + 4s\e^{\tt Z_{n}(V)}+
2\beta(s)\e^{\tt Z_{n}(V)} (4\tt\eta_{n,k} + \tt
\gg_k+4\ll\mu(\rr>n-k))\big\}.\end{split}\end{equation*} So, for any
$r>0$, let $s= \ff r 8 \e^{-\tt Z_{n}(V)}$. If for some $n,k\ge 1$
one has
$$ 6\mu_V(\rr>n)
   +2\e^{\tt Z_{n}(V)}\beta\Big(\ff r 8 \e^{-\tt Z_{n}(V)}\Big)\big(4\tt \eta_{n,k} +\tt \gg_k+ 4\ll\mu(\rr>n-k)\big)\le \ff r 2,$$ then
$$\mu_V(f^2)\le 2\beta(s)\e^{\tt K_{n,k}(V)} \EE_V(f)+r\|f\|_\infty^2.$$ Therefore, the proof is finished.
\end{proof}

To conclude this section,  we present an example where $(\EE,
\scr{D}(\EE))$ satisfies the Poincar\'e inequality, i.e. the weak
Poincar\'e inequality \eqref{wp} holds for a constant function
$\bb.$

\begin{exa}\label{exa-wp}\rm Let $E=\R^m$  with $d(x,y)=|x-y|$.
Let
$$q(x,y)=\ff{(1+|y|)^{m+\aa}}{|x-y|^{m+\alpha}},\ \  \mu(\d x)=\frac{c_{m,\aa}\d
x}{(1+|x|)^{m+\alpha}}$$ for some constant $0<\alpha<2$, where
$c_{m,\aa}$ is a normalizing constant such that $\mu$ is a
probability measure. Then $\A\supset C_0^2(\R^m),$ and according to \cite[Corollary 1.2(1)]{WW}, \eqref{wp} holds
for a constant rate function $\beta(r)\equiv \beta>0.$

Now, let   $V$ be measurable satisfying
\begin{equation}\label{exa-wp-01}
-s\varepsilon\log(1+|x|)-K \le  V(x)\le
(1-s)\varepsilon\log(1+|x|)+K,\quad x\in \R^m
\end{equation}
for some constants $\varepsilon\in [0,\alpha)$, $s \in [0,1]$ and
$K\in \R$ such that $\mu(e^V)=1$. Then the weak Poincar\'{e} inequality \eqref{th-wp-2} holds
with
\begin{equation}\label{exa-wp-02}
\beta_V(r)=C\big(1+r^{-\varepsilon/(\aa-(1-s)\vv)}\big)
\end{equation}
for some constant $C>0$.

Moreover, the assertion is sharp in the following two cases with $s=0$.
\beg{enumerate}\item[(i)] $\beta_V$ in     \eqref{exa-wp-02} is sharp, i.e. $\beta_V$ can not be replaced by any essentially smaller functions:    if
$$\lim_{r\to0} r^{\vv/(\aa-\vv)}\beta_V(r) =0,$$ then the weak Poincar\'{e} inequality \eqref{th-wp-2} does not hold.
 \item[(ii)] The constant $K$ can not be replaced by any unbounded functions: for \begin{equation}\label{exa-wp-03}
 V(x)=
 \varepsilon\log(1+|x|)+ \phi(|x|)+K_0,
\end{equation}
where $\vv\in [0,\aa)$,
$\phi:[0,+\infty)\rightarrow [0,+\infty)$
is an increasing function with $\phi(r)\uparrow\infty$ as $r\uparrow\infty$ such that
$\mu(e^{\varepsilon\log(1+|\cdot|)+ \phi(|\cdot|)})<\infty$, and $K_0\in\R^d$ is such that $\mu(e^V)=1$, the weak Poincar\'e inequality
(\ref{th-wp-2}) with the rate function $\beta_{V}$ given by
(\ref{exa-wp-02}) does not hold. \end{enumerate}
\end{exa}

\begin{proof} Take $k=\frac{n}{2}$. Then
  there is a constant $c_1>0$ such that for $n$ large enough
\begin{equation*}
\begin{split}
&\tt K_{n,\frac{n}{2}}(V)  \le \vv \log(1+n)+c_1,\ \ \ \tt Z_n(V)\le \vv(1-s) \log(1+n)+c_1,\\
&\tt\eta_{n,\frac{n}{2}}+\tt \gg_{\frac{n}{2}}  +\mu\big(\rr>\frac{n}{2}\big)\le \frac{c_1}{n^{\aa}},
\ \ \mu_V(\rr>n)\le \frac{ c_1}{n^{\aa-(1-s)\vv}}.
\end{split}
\end{equation*}
Since $\vv \in [0,\aa)$, we see that (\ref{th-wp-1}) holds and there exists $c_2>0$ such that
$$6\mu_V(\rr>n)
   +2\beta\e^{\tt Z_{n}(V)} \big(4\tt \eta_{n,\frac{n}{2}}+\tt \gg_{\frac{n}{2}} +4\ll\mu(\rr>\frac{n}{2})\big)\le \ff{c_2}{2n^{\aa-(1-s)\vv}},\ \ n\ge 1.$$
Thus, in the definition of $\beta_V$ for small $r>0$ we may take $n= (\ff{c_2}r)^{\ff 1{\aa- (1-s)\vv}}$ to get
 $$\beta_V(r)\le 2\beta \e^{\tt K_{n,\frac{n}{2}}(V)}\le c_3r^{-\vv/(\aa-(1-s)\vv)}$$ for some constant $c_3>0.$
Therefore, there exists $C>0$ such that the weak Poincar\'{e} inequality \eqref{th-wp-2} holds for   $\beta_V$ given in
\eqref{exa-wp-02}.

It remains to verify (i) and (ii), where the assertion (i) has been included in \cite[Corollary 1.2(3)]{WW}. So, it suffices to consider (ii).
  Now, assume
that the weak Poincar\'e inequality holds for
$(\EE_{V},\scr{D}(\EE_{V}))$ with the rate function
$$
\beta_V(r)=c_4(1+r^{-\frac{\vv}{\alpha-\vv}})
$$ for some constant $c_4>0$. For any $n\ge1$, let $f_n\in C^\infty(\R^d)$ be in part (b) of the proof of Example \ref{ex-sp-1}. By (\ref{AB0}) and noting that $\mu_V(f_n)^2\le \ff 1 2 \mu_V(f_n^2)$ for large $n$, there exist constants $c_5, c_6>0$ such that for $n$ large enough,
  $$\EE_V(f_n)\le c_5 n^{-\aa},\quad \mu_V(f_n^2)-\mu_V(f_n)^2\ge \frac{c_6\e^{\phi(n)}}{n^{\aa-\vv}}.$$ Combining these with \eqref{th-wp-2}, we obtain that there exists $c_7>0$ such that for all $r>0$ and for $n$ large enough,
  $$\frac{c_7\e^{\phi(n)}}{n^{\aa-\vv}}\le \frac{\bb_V(r)}{n^\aa}+r.$$
 Taking $r=\frac{c_7}{2}\frac{\e^{\phi(n)}}{n^{\aa-\vv}}$ in the inequality above, we arrive at
\begin{equation}\label{exa-wp-04}\bb_V\Big(\frac{c_7}{2}\frac{\e^{\phi(n)}}{n^{\aa-\vv}}\Big) \ge \frac{c_7}{2}n^\vv \e^{\phi(n)}.\end{equation}

Since there is $c_8>0$ such that for $n$ large enough
$$\frac{\e^{\phi(n)}}{n^{\aa-\vv}}\le c_8\int_{\{|x|\ge n\}}\frac{\e^{\phi(|x|)}}{(1+|x|)^{m+\aa-\vv}}
\,\d x,$$ it holds that $$\lim_{n\to \infty}\frac{\e^{\phi(n)}}{n^{\aa-\vv}}=0,$$ which, along with the definition of $\bb_V$, yields that there is a constant $c_9>0$ such that for $n$ large enough
$$ \bb_V\Big(\frac{c_7\e^{\phi(n)}}{2n^{\aa-\vv}}\Big)\le c_9n^\vv \e^{-\vv\phi(n)/(\aa-\vv)},$$
 which is a contradiction to \eqref{exa-wp-04} since $\lim_{r\to\infty}\phi(r)=\infty.$
 Therefore, the weak Poincar\'e
inequality does not hold with the rate function
(\ref{exa-wp-02}).
\end{proof}

\paragraph{Acknowledgement.}  We are grateful to an
anonymous referee whose comments helped to improve the exposition of this paper.

\end{document}